\documentclass[11pt]{article}
\usepackage{amssymb,latexsym}
\usepackage[dvips]{graphicx} 
\usepackage{url}

\newtheorem{prp}{Proposition}
\newtheorem{lem}[prp]{Lemma}\newtheorem{thm}[prp]{Theorem}
\newenvironment{prf}{\begin{trivlist}\item[\emph{Proof.}]}{\end{trivlist}
  \medskip\par}
\newenvironment{rem}{\begin{trivlist}\item[\emph{Remarks.}]}{\end{trivlist}
  \medskip\par}
\def\prpb{\begin{prp}}\def\prpe{\end{prp}}
\def\lemb{\begin{lem}}\def\leme{\end{lem}}
\def\thmb{\begin{thm}}\def\thme{\end{thm}}
\def\corb{\begin{cor}}\def\core{\end{cor}}
\def\prfb{\begin{prf}}\def\prfe{\end{prf}}
\def\remb{\begin{rem}}\def\reme{\end{rem}}
\def\prpa#1{\label{p:#1}}\def\prpu#1{Proposition~\ref{p:#1}}
\def\lema#1{\label{l:#1}}\def\lemu#1{Lemma~\ref{l:#1}}
\def\thma#1{\label{t:#1}}\def\thmu#1{Theorem~\ref{t:#1}}

\def\seca#1{\label{s:#1}}\def\secu#1{Section~\ref{s:#1}}

\def\itmb{\begin{enumerate}}\def\itme{\end{enumerate}}
\def\itdb{\begin{itemize}}\def\itde{\end{itemize}}
\def\ittb{\begin{description}}\def\itte{\end{description}}
\def\eqnb{\begin{equation}}\def\eqne{\end{equation}}
\def\arrb#1{\begin{array}{#1}}\def\arre{\end{array}}
\def\tabb#1{\par\noindent\begin{tabular}{#1}}
\def\tabe{\end{tabular}\par\noindent}
\def\eqna#1{\label{e:#1}}\def\eqnu#1{(\ref{e:#1})}

\def\QED{\relax\ifmmode\let\@tempa\relax\ifcase\@eqcnt\def\@tempa{& & &}\or
  \def\@tempa{& &}\else\def\@tempa{&}\fi\@tempa $\Box$ \else\hfill $\Box$ \fi}
\def\DDD{\relax\ifmmode\let\@tempa\relax\ifcase\@eqcnt\def\@tempa{& & &}\or
 \def\@tempa{& &}\else\def\@tempa{&}\fi\@tempa $\Diamond$
 \else\hfill $\Diamond$ \fi}

\def\Rom#1{\uppercase\expandafter{\romannumeral#1}}
\def\dsp{\displaystyle}
\def\eps{\epsilon}

\def\limf#1{\displaystyle \lim_{#1\to\infty}}

\def\Ccomb#1#2{\setbox0=\hbox{$\displaystyle\mathrm{C}$}\setbox1=\hbox{%
$\scriptstyle #1$}\kern \wd1{\mathrm{C}}_{\kern -1.05\wd0\kern -0.99\wd1{#1}
 \kern 1.15\wd0{#2}}}

\def\clvec#1#2#3{\def\clvecone{#3}\left(\arrb{c} \dsp #1\\ \dsp #2
 \ifx\clvecone\empty\else\\ \dsp #3\fi\arre\right)}

\def\diff#1#2{\dsp\frac{d\,#1}{d#2}}

\def\pderiv#1#2{\dsp\frac{\partial\,#1}{\partial#2}}

\def\le{\leqq} \def\leq{\leqq}\def\ge{\geqq} \def\geq{\geqq}

\def\reals{{\mathbb R}}
\def\preals{{\mathbb R_+}}
\def\integers{{\mathbb Z}}\def\pintegers{{\mathbb Z}_+}

\def\prb#1{\def\prbone{#1}
  \ifx\prbone\empty{\mathrm{P}}\else{\mathrm{P[\;}}#1{\mathrm{\;]}}\fi}
\def\prbseq#1#2{\def\prbseqone{#2}
  \ifx\prbseqone\empty{\mathrm{P}}_{#1}\ignorespaces
  \else{\mathrm{P}}_{#1}{\mathrm{[\;}}#2{\mathrm{\;]}}\fi}
\def\EE#1{{\mathrm{E[\;}}#1{\mathrm{\;]}}}
\def\EEseq#1#2{\def\EEseqone{#2}
  \ifx\EEseqone\empty{\mathrm{E}}_{#1}\else
 {\mathrm{E}}_{#1}{\dsp\mathrm{[\;}}#2{\mathrm{\;]}}\fi}

\def\VVseq#1#2{\def\VVseqone{#2}
  \ifx\VVseqone\empty{\matrm{V}}_{#1}\else
 {\mathrm{V}}_{#1}{\dsp\mathrm{[\;}}#2{\mathrm{\;]}}\fi}

\def\ssN{^{(N)}}

\title{
 Existence of an infinite particle limit of stochastic ranking process
}
\author{
Kumiko Hattori\thanks{
Department of Mathematics and Information Sciences,
 Tokyo Metropolitan University, Hachioji, Tokyo 192-0397, Japan.
\hspace*{1em} email: khattori@tmu.ac.jp
}
\\ \and
Tetsuya Hattori\thanks{
Mathematical Institute, Graduate School of Science, Tohoku University,
 Sendai 980-8578, Japan.
\hspace*{1em} URL: http://www.math.tohoku.ac.jp/\~{}hattori/amazone.htm
\hspace*{1em} email: hattori@math.tohoku.ac.jp
}
}
\date{\today}
\begin{document}
\maketitle

\begin{center}
ABSTRACT
\end{center}
We study a stochastic particle system which models the time evolution
of the ranking of books by online bookstores (e.g., Amazon).  In this system, 
particles are lined in a queue.  Each particle jumps at random jump 
times to the top of the queue, and otherwise stays in the queue, 
being pushed toward the tail every time another particle jumps to the top.
In an infinite particle limit, the random motion of each particle 
between its jumps converges to a deterministic trajectory.
(This trajectory is actually observed in the ranking data on web sites.)
We prove that the (random) empirical 
distribution of this particle system 
converges to a deterministic space-time dependent distribution.
A core of the proof is the law of large numbers for {\it dependent} 
random variables.

\vspace*{1in}\par
\bigskip\par
\noindent\textit{MSC2000 Subject Classifications:}
82C22, 60K35, 91B02
\bigskip\par
\noindent\textit{Corresponding author:} 
Tetsuya Hattori, 
hattori@math.tohoku.ac.jp
\par\noindent
Mathematical Institute, Graduate School of Science, Tohoku University,
 Sendai 980-8578, Japan
\par\noindent
tel+FAX:  011-81-22-795-6391

\newpage

\section{Introduction.}
\seca{1}

\subsection{Definitions.}

Let $(\Omega,{\cal B},\prb{})$ be a probability space, and
on this probability space we consider a stochastic ranking process
$\{X\ssN_i(t)\mid t\ge0,\ i=1,2,\cdots,N\}$ of $N$ particles,
where $X\ssN_i(t)$ is the ranking of particle $i$ at time $t$,
defined as follows.

For each $i$ assume that $x\ssN_{i,0}\in \{1,2,\cdots,N\}$ and $w\ssN_i>0$
are given.
$x\ssN_{i,0}$ is the initial value of the ranking of the particle $i$;
$x\ssN_{i,0}=X\ssN_i(0)$.
We require that $x\ssN_{i,0}$\,, $i=1,2,\cdots,N$, are
different numbers, or in other words,
$x\ssN_{i,0}$, $i=1,2,\cdots,N$, is a permutation of $1,2,\cdots,N$.
$w\ssN_i$ is the jump rate of the particle $i$.

For each $i$ let
$\tau\ssN_{i,j}$, $j=0,1,2,\cdots$,
be an increasing sequence of random jump times,
such that 
$\{\tau\ssN_{i,j}\mid j=0,1,2,\cdots\}$, $i=1,2,\cdots,N$, are independent
(independence among particles),
$\tau\ssN_{i,0}=0$ and
$\{\tau\ssN_{i,j+1}-\tau\ssN_{i,j}\mid j=0,1,2,\cdots \}$
are i.i.d.\ with the law of $\tau\ssN_i=\tau\ssN_{i,1}$ being
\eqnb
\eqna{stoppingtimeexponential}
\prb{\tau\ssN_i\le t} =1-e^{-w\ssN_it},\ t\ge 0\,.
\eqne
Note that with probability $1$,
$\tau\ssN_{i,j}$, $j=0,1,2,\cdots$, is strictly increasing,
and that 
$\tau\ssN_{i,j} \ne \tau\ssN_{i',j'}$
for any different pair of suffices $(i,j)\ne (i',j')$.

For each $i=1,2,\cdots,N$ we define the time evolution of $X\ssN_i$ by,
\eqnb
\eqna{defprocjx1}
X\ssN_i(t)=x\ssN_{i,0}+
\sharp\{i'\in\{1,2,\cdots,N\} \mid x\ssN_{i',0}>x\ssN_{i,0},
\ \tau\ssN_{i',1} \le t \},
\ \ 0\le t< \tau\ssN_{i,1},
\eqne
where $\sharp A$ denotes the number of elements in the set $A$,
with $\sharp\emptyset=0$,
and for each $j=1,2,3,\cdots$
\eqnb
\eqna{defprocjx}
\arrb{l}\dsp
X\ssN_i(\tau\ssN_{i,j})=1,\ \mbox{ and }
\\ \dsp
X\ssN_i(t)=
\sharp\{i'\in\{1,2,\cdots,N\} \mid
 \exists j'\in\pintegers;
\ \tau\ssN_{i,j} <  \tau\ssN_{i',j'} \le t \},
\ \ \tau\ssN_{i,j}<t<\tau\ssN_{i,j+1}.
\arre
\eqne

Intuitively speaking, the definition says that
particle $i$ jumps at random times $\tau_{i,j}$
to the top of the queue, and that after the jump it is
pushed toward the tail every time another particle 
of larger ranking number jumps to the top.
For example, let $N=4$ and let the initial ranking
be $x^{(4)}_{1,0}=2$, $x^{(4)}_{2,0}=3$, $x^{(4)}_{3,0}=1$, $x^{(4)}_{4,0}=4$.
In other words, particles $1$--$4$ are initially aligned as $3124$.
For a sample $\omega$ such that
\[ 0<\tau_{1,1}(\omega)<\tau_{2,1}(\omega)<
\tau_{4,1}(\omega)< \tau_{1,2}(\omega) < \cdots, \]
the configuration evolves as
\[ 3124\ \to\ 1324\ \to\ 2134\ \to\ 4213\ \to\ 1423\ \to\ \cdots, \]
where the changes occur at each jump times $\tau_{i,j}(\omega)$.

The stochastic ranking process may be viewed as a mathematical model of 
the time evolution of rankings such as that of books
on the online bookstores' web (e.g., www.Amazon.co.jp).
In this example, 
$N$ stands for the total number of books,
$i$ represents a specific title of a book,
$w\ssN_i$ is the average rate with which the book $i$ is sold,
$x\ssN_{i,0}$ is the initial position (ranking) of the book,
$\tau\ssN_{i,j}$ is the random time at which the book $i$ is
sold for the $j$-th time,
and $X\ssN_i(t)$ is the ranking of the book $i$ at time $t$.

In the time interval $(\tau\ssN_{i,j},\tau\ssN_{i,j+1})$
the ranking $X\ssN_i(t)$ increases by $1$ every time one of
the books in the tail side of the ranking 
(i.e., with larger $X\ssN_{i'}(t)$) 
is sold. In other words we have the following.
\prpb
\prpa{defprocd}
\eqnu{defprocjx1} and \eqnu{defprocjx} are equivalent to the following:
For each $i$,
\itmb
\item
$X\ssN_i(\tau\ssN_{i,j})=1$, $j=1,2,\cdots$,
\item
for each $i'\ne i$ and $j'=1,2,\cdots$,
if $X\ssN_i(\tau\ssN_{i',j'}-0)<X\ssN_{i'}(\tau\ssN_{i',j'}-0)$ then
$\dsp X\ssN_i(\tau\ssN_{i',j'})=X\ssN_i(\tau\ssN_{i',j'}-0)+1$,
\item
otherwise $X\ssN_i(t)$ is constant in $t$.
\DDD
\itme
\prpe

As seen in \prpu{defprocd}, each particle jumps at random times to rank $1$,
and gradually moves to the right (increasing number) 
without outpacing any other particles on its right.
This implies that for each $t$ 
there is a boundary position $x\ssN_C(t)\in\{0,1,\cdots,N-1\}$
such that all the particles on the left side have experienced a jump,
and that none of the particles on the right has jumped by time $t$:
\[
\arrb{l@{\ \Rightarrow\ }l} \dsp X\ssN_i(t)<x\ssN_C(t) & \tau\ssN_i\le t,
 \\ X\ssN_i(t)\ge x\ssN_C(t) & \tau\ssN_i>t. \arre
\]
$x\ssN_C(t)$ is a random variable and is explicitly written as: 
\eqnb
\eqna{xNCt}
x\ssN_C(t)= \sum_{i=1}^N \chi_{\tau\ssN_i\le t}\,.
\eqne
Put
\eqnb
\eqna{yNCt}
y\ssN_C(t)=\frac1N x\ssN_C(t)=\frac1N \sum_{i=1}^N \chi_{\tau\ssN_i\le t}\,.
\eqne

\subsection{Motivation.}

We are interested in the large $N$ limit of the stochastic ranking process.
Noting that \eqnu{stoppingtimeexponential} implies
$\dsp \EE{\chi_{\tau\ssN_i\le t}}=\prb{\tau\ssN_i\le t}=1-e^{-w\ssN_i t}$,
the weak law of large numbers easily leads to the following:
\prpb
\prpa{yCt}
Assume that the empirical distribution of jump rates
converges to a probability distribution $\lambda$:
\eqnb
\eqna{HDLlimitPareto}
\lambda\ssN(dw)=\frac1N \sum_{i=1}^N \delta(w-w\ssN_i)\,dw
 \to \lambda(dw), \ N\to\infty.
\eqne
Then
\[
\limf{N} y_C\ssN(t)=y_C(t),\ \ \mbox{ in probability,}
\]
where
\eqnb
\eqna{yCt}
y_C(t)=1-\int_{[0,\infty)} e^{-wt} \lambda(dw),\ \ t\ge0.
\eqne
\DDD\prpe
In the case of the online bookstore Amazon.com, the rankings 
seem to be defined in a more involved way, but
the trajectories of rankings as predicted by \prpu{yCt} are actually
observed at Amazon.co.jp \cite{HH072,HH073}.
As may be seen from this example,
the stochastic ranking process
would be of increasing practical interest and significance in this
age of online retails and web 2.0 \cite{longtail}.

In this paper we go further and prove that in the infinite particle limit
$N\to\infty$,
the random empirical distribution of the particle system converges
to a deterministic space-time dependent distribution.
To consider the limit $N\to\infty$, it is natural to use the
spacially scaled variables:
\eqnb
\eqna{defprocj}
Y\ssN_i(t)=\frac1N\, (X\ssN_i(t)-1).
\eqne
$Y\ssN_i(t)$ denotes the
spacially scaled position of the particle $i$ at time $t$, taking 
values in $\dsp [0,1)\cap N^{-1}\integers$.
In the following, we will use 
$Y\ssN_i(t)$ instead of $X\ssN_i(t)$.
Correspondingly, we will use
\[ y\ssN_{i,0}=\frac1N\, (x\ssN_{i,0}-1)\in [0,1)\cap N^{-1}\integers, \]
for the initial configurations instead of $x\ssN_{i,0}$.
\thmb[\protect\thmu{HDL}]
Under the assumptions on the initial configurations in \secu{21},
the joint random empirical distribution of jump rates (particle types)
and positions associated with the stochastic ranking process $\{Y\ssN_i\}$
converges as $N\to\infty$ to a distribution 
(with deterministic time evolution).
\DDD\thme
See \thmu{HDL} in \secu{classicalfieldrelaxationconserved} for
a precise form of the statement and an explicit form of
the limit distribution.

In the ranking of books,
each time a book is sold its ranking jumps to $1$,
no matter how unpopular the book may be.
At first thought one might guess that such a naive ranking
will not be a good index for the popularity of books.
But thinking more carefully, one notices that the well sold books
(particles with large $w\ssN_i$, in our definition) are dominant
near the top position, while books near the tail position
are rarely sold.
Hence, though the rankings of each book is stochastic and has sudden jumps,
the spacial \textit{distribution} of jump rates are more stable,
with the ratio of books with large jump rate high near the top position
and low near the tail position.
Seen from the bookstore's side, it is not a specific book that really
matters, but a totality of book sales that counts, so 
the evolution of distribution of jump rate should be important.
An intuitive meaning of the Theorem is that we can make
this intuition rigorous and precise, with an explicit form of the 
distribution for large $N$ limit.

The limit in the Theorem is mathematically non-trivial in that it involves
the law of large numbers for dependent variables.
Dependence occur because,
for each particle $i$, 
the time evolution between the jump times $\tau\ssN_i$ is
a trajectory of a flow caused by the jumps of other particles
in the tail side of the ranking,
and the conditioning on tail side induces dependence
of stochastic variables.

The idea of considering such a limit theorem is mathematically motivated
by the celebrated theory of hydrodynamic limits \cite{KL,OVY93,Spohn,VY97},
although the dynamics (relaxation to equilibrium) and hence the proofs in 
\secu{3} apparently have little in common (and are simpler).
A difference lies in  that the theory of hydrodynamic limits
(among other things) evaluate the relaxation to
equilibrium (invariant measures) through entropy and large deviation
arguments via local equilibrium, while
the dynamics of the stochastic ranking process has a special feature that
the queue of the particles consists of the `tail' regime
and the `head' regime, such that the former
is the queue of books which has not been sold up to time $t$, and
having no dynamics for relaxation,
keeps the remnant of initial data (\secu{ygttCtprf}).
In contrast, in the `head' regime
the `stationarity' is reached from the beginning
(\secu{ylttCtprf}).
It may also be worthwhile to note that our limit distributions,
unlike the hydrodynamic limits which satisfy diffusion equations,
satisfy non-local field equations (see remarks to \thmu{HDL} in \secu{2}).

In \secu{2} we state our main theorem and in \secu{3} we give a proof.

\smallskip\par\textbf{Acknowledgements.}
We thank the referee and the associate editor and Prof.~K.~Ishige
for constructive discussions, 
which (among other things) helped refining 
the assumptions for the main results.

The research of K.~Hattori 
is supported in part by a Grant-in-Aid for 
Scientific Research (C) 16540101 from the Ministry of Education,
 Culture, Sports, Science and Technology, and
the research of T.~Hattori 
is supported in part by a Grant-in-Aid for 
Scientific Research (B) 17340022 from the Ministry of Education,
 Culture, Sports, Science and Technology.

\section{Main result.}
\seca{2}

\subsection{Assumptions on initial configuration.}
\seca{21}

We consider the $N\to\infty$ limit of the empirical distribution
on the product space of jump rate and spacial position $\preals\times[0,1)$,
\eqnb
\eqna{HDLtN}
\mu\ssN_{t}(dw,dy)
=\frac1N \sum_i \delta(w-w\ssN_i)\delta(y-Y\ssN_i(t))\,dw\,dy.
\eqne

We impose that the initial distribution
\eqnb
\eqna{HDL0N}
\mu\ssN_{0}(dw,dy)
=\frac1N \sum_i \delta(w-w\ssN_i)\delta(y-y\ssN_{i,0})\,dw\,dy
\eqne
converges as $N\to\infty$ uniformly in $y$ to 
a probability distribution
which is absolutely continuous with respect to the Lebesgue measure
in $y$.
(For the set of probability measures on $\preals\times[0,1)$
we work with the topology of weak convergence.)
Explicitly, we assume that for each $0\le y< 1$ there exists
a probability measure $\mu_{y,0}(dw)$ on $\preals$ such that
for any bounded continuous function $g:\ \preals\to\reals$,
the quantity
$\dsp \int_0^{\infty} g(w) \mu_{y,0}(dw)$ is Lebesgue measurable in $y$ and
\eqnb
\eqna{HDLinitialuniform}
\limf{N} \sup_{y\in [0,1)} \left|
\frac1N \sum_{i=1}^N g(w\ssN_i) \chi_{y\ssN_{i,0}\le y}
-\int_0^y \left(\int_0^{\infty} g(w) \mu_{z,0}(dw)\right)\,dz
\right|
=0,
\eqne
where (and also in the following) we use a notation $\chi_A$ which
is $1$ if $A$ is true and $0$ if $A$ is false.
That each $\mu_{y,0}(dw)$ is a probability measure on $\preals$ is
consistent with \eqnu{HDLinitialuniform}, which may be seen by 
letting $g(w)=1$, $w\ge 0$, in \eqnu{HDLinitialuniform} to find
\eqnb
\eqna{muy0normalization}
 \mu_{y,0}(\preals)=1,\ \ 0\le y<1\,.
\eqne
Note that
$\lambda\ssN$ and $\lambda$ in \eqnu{HDLlimitPareto}
are the marginal distribution of $\mu\ssN_{y,0}$ and $\mu_{y,0}$
of the jump rate;
\eqnb
\eqna{HDLinitialvslimitPareto}
\arrb{l}\dsp
\lambda\ssN(dw)
=\frac1N \sum_{i=1}^N \delta(w-w\ssN_i)\,dw
= \mu\ssN_{0}(dw , [0,1)),
\\ \dsp
\lambda(dw)=\int_0^1 \mu_{y,0}(dw)\,dy.
\arre
\eqne
Note also that \eqnu{HDLinitialuniform} and Fubini's Theorem imply
\eqnu{HDLlimitPareto}.

We assume that the average of $\lambda$ is finite,
\eqnb
\eqna{lambdafiniteaverage}
\int w \lambda(dw)<\infty,
\eqne
and 
\eqnb
\eqna{lambda00}
\lambda(\{0\})=0.
\eqne
This complets the assumptions for our main results.

\remb
\itmb
\item
The assumption \eqnu{lambdafiniteaverage}
assures that $\mu_{0,t}$ is well-defined (see \eqnu{Tets20070726}).
The main results on the existence of the infinite particle limit will hold 
without \eqnu{lambdafiniteaverage} for $y>0$,
but we keep this assumption to include $y=0$.

\item
We assume \eqnu{lambda00} to assure that $y_C:\ [0,\infty)\to[0,1)$ 
defined in \eqnu{yCt} is onto (see \prpu{relaxationcurve}).
The basic results in this paper will hold (with extra complexity in
notations and arguments) without \eqnu{lambda00},
but we prefer to keep notations and arguments simple by keeping this
assumption.
\eqnu{lambda00} implies in the actual bookstore ranking, 
that `all the books sell (almost surely)'.
\DDD
\itme
\reme

\subsection{Main theorem.}
\seca{classicalfieldrelaxationconserved}

With \eqnu{lambda00}, it is straightforward to show  
\prpb
\prpa{relaxationcurve}
Assume \eqnu{lambda00}. Then
$y_C:\ [0,\infty)\to[0,1)$ defined in \eqnu{yCt}
is a continuous, strictly increasing, bijective 
function of $t$.
\DDD\prpe
\prpu{relaxationcurve} implies the existence of the inverse function  
$t_0:\ [0,1)\to[0,\infty)$, satisfying
\eqnb
\eqna{t0y}
y_C(t_0(y))=y, \ \ 0\le y<1\,,
\eqne
or 
\eqnb
\eqna{yCasevaporatedparticlest0}
y=1-\int_0^{\infty} e^{-wt_0(y)} \lambda(dw).
\eqne

Differentiating \eqnu{yCt} and \eqnu{t0y}, we have
\eqnb
\eqna{dyCdt}
\diff{y_C}{t}(t)= \int_0^{\infty} we^{-wt} \lambda(dw)
=\frac{1}{\dsp \diff{t_0}{y}(y_C(t))}  \,.
\eqne

We generalize \eqnu{yCt} and define (with slight abuse of notations)
\eqnb
\eqna{ytildeCyt}
y_C(y,t)=1-\int_y^1 \int_0^{\infty} e^{-wt} \mu_{z,0}(dw)\,dz,
\ \ t\ge0,\ 0\le y< 1\,.
\eqne
In particular, $y_C(t)={y}_C(0,t)$.
In the infinite  particle limit, 
${y}_C(y,t)$ denotes the position of a particle at time $t$
(if it does not jump up to time $t$)
 whose initial position is $y$ (\prpu{trajectory}).
 
\prpb
\prpa{ycytstrictincrease}
${y}_C(\cdot,t):\ [0,1)\to[y_C(t),1)$ is a continuous, 
strictly increasing,
bijective function of $y$. 
\DDD\prpe
\prfb
It is straightforward from the definition of $y_C(y,t)$ in \eqnu{ytildeCyt}
to see 
that ${y}_C(\cdot,t)$ is continuous and non-decreasing in $y$.
To see that it is strictly increasing, let $0\leq z_2 <z_1<1$ .
Then \eqnu{ytildeCyt} implies
$\dsp
y_C(z_1,t)-y_C(z_2,t)
=\int_{z_2}^{z_1} \int_0^{\infty} e^{-wt} \mu_{z,0}(dw)\,dz.
$
If this is $0$, 
then $\mu_{z,0}([0,M])=0$ for any $M>0$, for a.e.~$z\in [z_2,z_1]$,
which contradicts \eqnu{muy0normalization}.
\QED\prfe
  
\prpu{ycytstrictincrease} implies that 
the inverse function 
$\hat{y}(\cdot,t):\ [y_C(t),1)\to[0,1)$ exists:
\eqnb
\eqna{yhatyt}
1-y=\int_{\hat{y}(y,t)}^1 \int_0^{\infty} e^{-wt} \mu_{z,0}(dw)\,dz,
\ \ t\ge0,\ y_C(t)\le y< 1.
\eqne
$\hat{y}(y,t)$ denotes the initial position of a particle
located at $y$ ($>y_C(t)$) at time $t$.
It holds that

\eqnb
\eqna{hatyy}
\pderiv{\hat{y}}{y}(y,t)=\frac{1}{\dsp
 \int_0^{\infty} e^{-wt} \mu_{\hat{y}(y,t),0}(dw)}\,.
\eqne

Now we return to our $N$-particle process.
\thmb
\thma{HDL}
Consider  the stochastic ranking process
$\{Y\ssN_i\}$ defined by \eqnu{stoppingtimeexponential} and
\eqnu{defprocj}.
Assume 
\eqnu{HDLinitialuniform}
 \eqnu{lambdafiniteaverage} and \eqnu{lambda00}.
Then the joint empirical distribution of particle types and positions 
at time $t$ 
\eqnb
\eqna{empiricaldistribution}
\mu\ssN_{t}(dw,dy)=
\frac1N \sum_i \delta(w-w\ssN_i)\delta(y-Y\ssN_i(t))\,dw\,dy
\eqne
converges as $N\to\infty$ to a distribution 
$\mu_{y,t}(dw)\,dy$ on $\preals\times [0,1)$;
for any bounded continuous function 
$f:\ \preals\times[0,1) \to \reals$
\eqnb
\eqna{HDL}
\arrb{l}\dsp
\limf{N} \frac1N \sum_i f(w\ssN_i,Y\ssN_i(y) )
= \int_0^1 \left(\int_0^{\infty} f(w,y) \mu_{y,t}(dw)\right)\,dy,
\ \ \mbox{ in probability. }
\arre \eqne
 
The limit distribution
is absolutely continuous with respect to the Lebesgue measure on $[0,1)$.
The density $\mu_{y,t}(dw)$ with regard to $y$ is given by
\eqnb
\eqna{Tets20070726}
\mu_{y,t}(dw)=\left\{ \arrb{ll}\dsp
\frac{\dsp
we^{-wt_0(y)} \lambda(dw)
}{\dsp
\int_0^{\infty} we^{-wt_0(y)} \lambda(dw)
}\,, & y<y_C(t), \\
\frac{\dsp
e^{-wt} \mu_{\hat{y}(y,t),0}(dw)
}{\dsp
\int_0^{\infty} e^{-wt} \mu_{\hat{y}(y,t),0}(dw)
}\,, & y>y_C(t).
\arre \right. 
\eqne
\DDD\thme

\remb
\itmb
\item

\eqnu{empiricaldistribution} and \eqnu{HDLinitialvslimitPareto}
imply
$\dsp \lambda\ssN(\cdot)= \mu\ssN_{t}(\cdot,[0,1))$.
Moreover,
if in \eqnu{HDL} we take $f$ without $y$ dependence
and use \eqnu{HDLlimitPareto}, 
we have as a generalization of 
\eqnu{HDLinitialvslimitPareto}
\eqnb
\eqna{HDLvslimitPareto}
\lambda=\int_0^1 \mu_{y,t}\,dy .
\eqne

\item
Our results state that a random phenomenon approaches a deterministic one
as the particle number $N$ is increased. We state the results in
terms of convergence in probability, but since the limit quantity
 is deterministic, this limit is equivalent to convergence in law.

\item
The explicit forms in \eqnu{Tets20070726} differ drastically for
$y>y_C(t)$ and $y<y_C(t)$. As we have pointed out in the Introduction,
and also as we will see in the proofs in the next section,
the dynamics for large $y$ and small $y$ are totally different.

\item
By direct calculations,
 one sees that $\mu_{y,t}(dw)$ satisfies the following (non-local) equations:
\eqnb
\eqna{fluid}
\pderiv{\mu_{y,t}(dw)}{t} +\pderiv{(v(y,t)\, \mu_{y,t}(dw))}{y}
=-w \mu_{y,t}(dw),
\eqne
where
\eqnb
\eqna{velocity}
 \arrb{l}\dsp
v(y,t)
=\int_{y}^1 \left( \int w\, \mu_{z,t}(dw) \right)\,dz
=\left\{ \arrb{ll} \dsp
\pderiv{y_C}{t}(t_0(y))\,,
& y<y_C(t), \\ \dsp
\pderiv{y_C}{t}(\hat{y}(y,t),t),
& y>y_C(t).
\arre \right.
\arre
\eqne
The partial differential equation \eqnu{fluid} can be seen 
as the equations of continuity (conservation of mass)
for the one-dimensional incompressible mixed fluids, 
with $w$ standing for the rate of evaporation of specific type of fluid
in the mixture \cite{HH072}.
$v(y,t)$ is the velocity of the fluid at position $y$ and time $t$,
and \eqnu{velocity}, the source of non-locality,
 means that the flow is driven by evaporation.
Intuitively, the equations are natural classical limit of
the stochastic processes considered in this paper.

In fact, we can directly prove (without referring to stochastic processes)
that \eqnu{Tets20070726} is the unique classical solution to the 
Cauchy problem of partial differential equation \eqnu{fluid} with
suitable boundary conditions.
See \cite{HH072} for details.
\DDD
\itme
\reme

\section{{Proof of \protect\thmu{HDL}.}}
\seca{3}

It is sufficient to consider the case that   
$f:\ \preals\times[0,1) \to \reals$ in  \eqnu{HDL} is 
expressed as $\dsp f(w,z)=g(w)\chi_{z\in [0,y]}$, 
with a bounded continuous function $g:\ \preals\to\reals$ 
and $0<y<1$. Thus we prove in this section
 \eqnb
\eqna{HDLequiv}
\limf{N} \frac1N \sum_i g(w\ssN_i) \chi_{Y\ssN_i(t)\le y}
= \int_0^y dz\, \int_0^{\infty} g(w) \mu_{z,t}(dw),\ \ \mbox{ in probability,}
\eqne
for any bounded continuous function $g:\ \preals\to\reals$.

\subsection{Case  ^^ $y=y_C(t)$'.}

\lemb
\lema{HDLyC}
For each $t>0$ and each bounded continuous function
$g\ :\ \preals\to\reals$, it holds that
\eqnb
\eqna{HDLyC}
\limf{N} \frac1N \sum_i g(w\ssN_i) \chi_{Y\ssN_i(t)\le y_C^{(N)}(t)}
=  \int_0^{\infty} g(w)(1- e^{-wt})\lambda(dw) , \ \ \mbox{ in probability,}
\eqne
and
\eqnb
\eqna{HDLyc2}
\limf{N} \frac1N \sum_i g(w\ssN_i) \chi_{Y\ssN_i(t)\le y_C^{(N)}(t)}
=\limf{N} \frac1N \sum_i g(w\ssN_i) \chi_{Y\ssN_i(t)\le y_C(t)},
\ \ \mbox{ in probability.}
\eqne
\DDD\leme
\prfb
Since we have from \eqnu{yNCt} and the definition of $Y\ssN_i(t)$
\eqnb
\eqna{HDLyCt0yprf1}
\chi_{Y\ssN_i(t)\le y\ssN_C(t)}
=
\chi_{\tau\ssN_i\le t},
\eqne
we have from \eqnu{stoppingtimeexponential}
\[
\EE{g(w\ssN_i) \chi_{Y\ssN_i(t)\le y\ssN_C(t)}}
= g(w\ssN_i)(1-e^{-w\ssN_it}).
\]
Since $g$ is bounded, in a similar manner as the proof of \prpu{yCt},
we apply the weak law of large numbers to obtain
\[
\limf{N} \left( \frac1N \sum_i g(w\ssN_i) \chi_{Y\ssN_i(t)\le y\ssN_C(t)}
-\frac1N \sum_i g(w\ssN_i)(1-e^{-w\ssN_it}) \right)
=0,\ \ \mbox{ in probability.}
\]
With \eqnu{HDLlimitPareto} we have \eqnu{HDLyC}.
Next, \prpu{yCt} implies that,
for any $\eps>0$, with large enough $N$, it holds that
\[ \prb{|y_C(t)- y_C\ssN(t)|\le\eps}> 1-\eps.\]
For such $N$, 
noting that 
\[
\sum_{i} | \chi_{Y\ssN_i(t)\le y\ssN_C(t)}-\chi_{Y\ssN_i(t)\le y_C(t)}  |
=
\sum_{k=1}^N  | \chi_{k\le Ny\ssN_C(t)}-\chi_{k\le Ny_C(t)}  |
\le N|y\ssN_C(t)-y_C(t)|+1,
\]
the probability of the event
\[\arrb{l}\dsp
 \left| \frac1N \sum_i g(w\ssN_i) \chi_{Y\ssN_i(t)\le y\ssN_C(t)}
-
\frac1N \sum_i g(w\ssN_i) \chi_{Y\ssN_i(t)\le y_C(t)} \right|
\\ \dsp
\le \frac1N \sup_w |g(w)| \sum_i
\left|  \chi_{Y\ssN_i(t)\le y\ssN_C(t)} -
 \chi_{Y\ssN_i(t)\le y_C(t)} \right|
\le \sup_w |g(w)| (\eps + \frac{1}{N})
\arre\]
is larger than $1-\eps$, thus 
the left-hand side converges to $0$ in probability as $N\to \infty$,
which implies \eqnu{HDLyc2}.
\QED\prfe

\subsection{Case $y<y_C(t)$.}
\seca{ylttCtprf}

First note that,
to prove \eqnu{HDLequiv} for $y<y_C(t)$, it is sufficient to prove that
for each bounded continuous function $g:\ \preals\to\reals$,
\eqnb
\eqna{HDLysmall}
\limf{N} \frac1N \sum_i g(w\ssN_i) \chi_{Y\ssN_i(t)\le y}
=  \int_0^{\infty} g(w) (1-e^{-wt_0(y)})\lambda(dw) ,
 \ \ \mbox{ in probability,}
\eqne
where $t_0(y)$ is as in \eqnu{t0y}.
To see that \eqnu{HDLysmall} implies \eqnu{HDLequiv},
differentiate the right-hand side of \eqnu{HDLysmall} with respect to $y$,
use \eqnu{dyCdt} and \eqnu{t0y}, and integrate from 
$0$ to $y$, keeping in mind  $t_0(0)=0$, and finally rewrite using
$\mu_{y,t}(dw)$ in \eqnu{Tets20070726},
to obtain
\[ \arrb{l}\dsp
\limf{N} \frac1N \sum_i g(w\ssN_i) \chi_{Y\ssN_i(t)\le y}
=  \int_0^y \frac{\dsp
\int_0^{\infty} g(w) we^{-wt_0(z)} \lambda(dw)
}{\dsp
\int_0^{\infty} we^{-wt_0(z)} \lambda(dw)
} \,dz
\\ \dsp
=  \int_0^y \int g(w)\mu_{z,t}(dw)\,dz,
 \ \ \mbox{ in probability,}
\arre \]
which gives \eqnu{HDLequiv}.

To prove \eqnu{HDLysmall},
fix $y<y_C(t)$ and let $t_0=t_0(y)<t$. 

Denote by $\{\tilde{Y}\ssN_i(s),\tilde{\tau}\ssN_{i,j}\}$ 
the scaled stochastic ranking process 
with the time origin shifted by the amount $t-t_0>0$\,.
Namely, let  $\tilde{Y}\ssN_i (s) = Y\ssN_i(s+t-t_0)$.
In particular, we have
\eqnb
\eqna{HDLprf11}
\frac1N \sum_i g(w\ssN_i) \chi_{Y\ssN_i(t)\le y}
=
\frac1N \sum_i g(w\ssN_i) \chi_{\tilde{Y}\ssN_i(t_0)\le y}\,.
\eqne
For $j=0,1,2,\cdots$, define $\tilde{\tau}\ssN_{i.j}$ by
$\tilde{\tau}\ssN_{i,0}=0$ and
\eqnb
\eqna{tildetau}
\tilde{\tau}\ssN_{i.j}=\tau\ssN_{i,\tilde{j}(i,t-t_0)+j-1} -(t-t_0)\,,
\eqne
where
\eqnb
\eqna{tildej}
\tilde{j}(i,t-t_0)=\inf\{ j\mid \tau\ssN_{i,j}> t-t_0\}.
\eqne
Put, in analogy to \eqnu{yNCt},
\[
\tilde{y}\ssN_C(s)=\frac1N \sum_{i=1}^N \chi_{\tilde{\tau}\ssN_i\le s}\,.
\]
Note that $\{\tau\ssN_{i,j+1}-\tau\ssN_{i,j}\mid j=0,1,2,\cdots\}$
are independent and have exponential distributions.
The loss of memory property of exponential distributions then implies
that $\{\tilde{\tau}\ssN_{i,j}\}$ have
 the same distributions as $\{\tau\ssN_{i,j}\}$,
and $\{\tilde{Y}\ssN_i(s) \}$ is a scaled stochastic ranking process
with jump times $\{\tilde{\tau}\ssN_{i,j}\}$ and initial configuration
$\{ Y\ssN_i(t-t_0) \}$.

Since $y=y_C(t_0)$,
 \eqnu{HDLyc2} for the time shifted ranking process implies
\eqnb
\eqna{HDLprf12}
\arrb{l}\dsp
\limf{N} \frac1N \sum_i g(w\ssN_i) \chi_{\tilde{Y}\ssN_i(t_0)\le y}
=\limf{N} \frac1N \sum_i g(w\ssN_i) \chi_{\tilde{Y}\ssN_i(t_0)\le y_C(t_0)}
\\ \dsp {}
=\limf{N} \frac1N \sum_i g(w\ssN_i) \chi_{\tilde{Y}\ssN_i(t_0)
\le \tilde{y}_C^{(N)}(t_0)},
\ \ \mbox{ in probability.}
\arre
\eqne
Using \eqnu{HDLyCt0yprf1}
 for the original process and the time shifted process, and recalling that
$\{\tau\ssN_i\}$ and $\{\tilde{\tau}\ssN_i\}$ have the same distribution,
and then using \eqnu{HDLyC}, we arrive at
\eqnb
\eqna{HDLprf13}
\arrb{l}\dsp
\limf{N} \frac1N \sum_i g(w\ssN_i) \chi_{\tilde{Y}\ssN_i(t_0)
\le \tilde{y}_C^{(N)}(t_0)}
=
\limf{N} \frac1N \sum_i g(w\ssN_i) \chi_{\tilde{\tau}\ssN_i\le t_0}
\\ \dsp
=\limf{N} \frac1N \sum_i g(w\ssN_i) \chi_{\tau\ssN_i\le t_0}
=
\limf{N} \frac1N \sum_i g(w\ssN_i) \chi_{Y\ssN_i(t_0)\le y_C^{(N)}(t_0)}
\\ \dsp {}
=  \int_0^{\infty} g(w)(1- e^{-wt_0})\lambda(dw) , \ \ \mbox{ in probability.}
\arre
\eqne

\eqnu{HDLprf11} \eqnu{HDLprf12} and \eqnu{HDLprf13} together imply
\eqnu{HDLysmall}.

\subsection{Case $y>y_C(t)$.}
\seca{ygttCtprf}

First we make some preparations for the proof.
We generalize $y\ssN_C(t)$ in \eqnu{yNCt} and define
(with slight abuse of notations) for $0\le y<1$
\eqnb
\eqna{HDLprf6}
y\ssN_C(y,t)=y+ \frac1N \sum_{i=1}^N \chi_{\tau\ssN_i\le t}
\chi_{y\ssN_{i,0}\ge y}.
\eqne
\prpb
\prpa{trajectory}
For $0\le y<1$ and $t\ge 0$,
\eqnb
\eqna{trajectory}
\limf{N} y\ssN_C(y,t) = y_C(y,t),\ \ \mbox{ in probability.}
\eqne
Namely, the (random) position $y\ssN_C(y,t)$ of a particle 
at time $t$ whose initial position is $y$ converges in probability
to a deterministic trajectory ${y}_C(y,t)$ defined by \eqnu{ytildeCyt}
in the infinite particle limit.
\DDD\prpe
\begin{trivlist}\item[\emph{Remark.}]
The proof below shows that the convergence in \eqnu{trajectory} holds 
uniformly in $y$.
\DDD\end{trivlist}\medskip\par
\prfb
\eqnu{HDLprf6} and \eqnu{stoppingtimeexponential} and
the independence of $\{\tau\ssN_i\}$ imply
\[ \arrb{l}\dsp
\EE{\left( y_C(y,t)-y\ssN_C(y,t)\right)^2}
\\ \dsp
=
 y_C(y,t)^2 
- 2y_C(y,t)\left(y+ \frac1N \sum_{i=1}^N 
(1-e^{-w\ssN_it})\chi_{y\ssN_{i,0}\ge y}\right)
\\ \dsp {}
+
y^2 + 2y\frac1N \sum_{i=1}^N (1-e^{-w\ssN_it})
\chi_{y\ssN_{i,0}\ge y}
\\ \dsp {}
+
 \frac1{N^2} \sum_{i=1}^N \sum_{j=1}^N
(1-e^{-w\ssN_it})(1-e^{-w\ssN_jt})
\chi_{y\ssN_{i,0}\ge y}
\chi_{y\ssN_{j,0}\ge y}
\\ \dsp {}
+
 \frac1{N^2} \sum_{i=1}^N
\left(1-e^{-w\ssN_it}-     (1-e^{-w\ssN_it})^2\right)
\chi_{y\ssN_{i,0}\ge y}
\\ \dsp
=
\left( y+ \frac1N \sum_{i=1}^N 
(1-e^{-w\ssN_it})\chi_{y\ssN_{i,0}\ge y} - y_C(y,t) \right)^2
\\ \dsp {}
+
 \frac1{N^2} \sum_{i=1}^N
\left(1-e^{-w\ssN_it}-     (1-e^{-w\ssN_it})^2\right)
 \chi_{y\ssN_{i,0}\ge y}\,.
\arre\]
The second term in the right-hand side of the equation above 
vanishes in the $N\to\infty$ limit because of the
factor $N^2$ in the denominator.
Concerning the first term, as
in the proof of \prpu{yCt} and \lemu{HDLyC}, 
\eqnu{HDLinitialuniform} implies that 
\[\arrb{l}\dsp
\limf{N}
\frac1N \sum_{i=1}^N 
(1-e^{-w\ssN_it})\chi_{y\ssN_{i,0}\ge y}
=\int_{y}^1 \int_0^{\infty} (1-e^{-wt}) \mu_{z,0}(dw)\, dz
\\ \dsp
=1-y- \int_{y}^1 \int_0^{\infty} e^{-wt} \mu_{z,0}(dw)\, dz,
\arre\]
uniformly in $y$. 
This combined with \eqnu{ytildeCyt} implies that the first term 
also vanishes.  Thus we have
\[\lim_{N\to \infty }\EE{\left( y_C(y,t)-y\ssN_C(y,t)\right)^2}=0.\]
With Chebyshev's inequality follows \eqnu{trajectory}.
\QED\prfe

As an equivalent statement to \eqnu{HDLequiv}, we will prove
\eqnb
\eqna{HDLequivl}
\limf{N} \frac1N \sum_i g(w\ssN_i) \chi_{Y\ssN_i(t)> y}
= \int_y^1 dz\, \int_0^{\infty} g(w) \mu_{z,t}(dw),
\ \ \mbox{ in probability.}
\eqne

Let  $\Omega\ssN_1=\{ y>y\ssN_C(t) \}$.  Then \prpu{yCt} implies that
if $y >y_C(t)$, 
\eqnb
\eqna{omega0}
\lim_{N \to \infty }\prb{\Omega\ssN_1} =1.
\eqne
For $t>0$ and $y>y\ssN_C(t)$, let
\eqnb
\eqna{yhatNyt}
\hat{y}\ssN(y,t)=\inf\{y\ssN_{i,0}\mid i=1,\cdots,N,\ Y\ssN_i(t)>y\}.
\eqne
Note that
\eqnb
\eqna{yNChatyNC}
 |y-y\ssN_C(\hat{y}^{(N)}(y,t),t)|\le \frac2{N}\,.
\eqne
This follows because
\eqnb
\eqna{yCN2YN}
 Y\ssN_i(t)=y\ssN_C(y\ssN_{i,0},t),\ \ \mbox{if } Y\ssN_i(t)>y\ssN_C(t),
\eqne
hence with $y>y\ssN_C(t)$
\[
y-\frac1N \le y\ssN_C(\hat{y}\ssN(y,t),t) \le Y\ssN_i(t)+\frac1N
\]
for all $i$ with $Y\ssN_i(t)>y$.

Until a particle jumps to the top of the queue, changes of its position
are caused only by the jumps of other particles that sit on its right
(\prpu{defprocd}), hence  
\eqnb
\eqna{HDLprf1}
 \sum_i g(w\ssN_i) \chi_{Y\ssN_i(t)> y}
= \sum_i g(w\ssN_i)
\chi_{ y\ssN_{i,0} \ge\hat{y}\ssN(y,t)}
\chi_{ \tau\ssN_i >t}, \ \ \mbox{ on } \Omega\ssN_1.
\eqne
Note that $\hat{y}\ssN(y,t)$ depends on  $\tau_i$'s. 
This means that the summands on the right-hand side are not
independent random variables, and that we can not apply the
law of large numbers as it is.
In contrast, since $\hat{y}(y,t)$ is deterministic, 
the law of large numbers yields,
as in the proofs of \prpu{yCt} and \lemu{HDLyC},
\eqnb
\eqna{HDLprf3}
\limf{N} \frac1N 
 \sum_i g(w\ssN_i)
\chi_{y\ssN_{i,0}\ge\hat{y}(y,t)}
\chi_{\tau\ssN_i>t}
= \int_{\hat{y}(y,t)}^1\, \int_0^{\infty} g(w) e^{-wt}\mu_{z,0}(dw)\, dz ,
\ \ \mbox{ in probability.}
\eqne
The right-hand side coincides with that of 
\eqnu{HDLequivl} 
through a change of variables $\hat{y}(z,t)\to z$
(note \eqnu{hatyy}).
Combining 
\eqnu{HDLequivl}, \eqnu{omega0}, \eqnu{HDLprf1} and 
\eqnu{HDLprf3}, 
we see that it is sufficient to prove that 
\eqnb
\eqna{HDLprf4}
\limf{N}
P [\ y>y\ssN_C(t), \ \ \frac1N 
 \sum_i g(w\ssN_i)
\left |\chi_{y\ssN_{i,0}\ge\hat{y}(y,t)}
-
\chi_{y\ssN_{i,0}\ge\hat{y}\ssN(y,t)} \right|
\chi_{\tau\ssN_i >t}\ge  \eps \ ]  =0\,
\eqne
holds for any $\eps \in(0,1)$.

Note that 
if $y>y\ssN_C(t)$ and 
$\dsp |\hat{y}\ssN(y,t)-\hat{y}(y,t)|\le \eps$,
then
\[\arrb{l}\dsp
\frac1N 
 \sum_i g(w\ssN_i)
\left |\chi_{y\ssN_{i,0}\ge\hat{y}(y,t)}
-
\chi_{y\ssN_{i,0}\ge\hat{y}\ssN(y,t)} \right|
\chi_{\tau_i>t}
\le
\frac{M}{N} \sum_i
\left |\chi_{y\ssN_{i,0}\ge\hat{y}(y,t)}
-
\chi_{y\ssN_{i,0}\ge\hat{y}\ssN(y,t)} \right|
\\ \dsp
\le M\,(\eps+\frac1N) ,
\arre \]
where $M>0$ is a constant satisfying $|g(w)|\le M$ for all $w\ge 0$\,.
Thus to prove
\eqnu{HDLprf4}, it is sufficient to show 
\[
\limf{N} \prb{y>y\ssN_C(t),\ \ |\hat{y}\ssN(y,t)-\hat{y}(y,t)|> \eps}=0,
\]
for an arbitrary $\eps >0$.

Let us assume otherwise, that is, with \eqnu{omega0} in mind, assume that  
there is $\eps_1 >0$ ,  $\rho >0$  and an increasing sequence of 
positive integers $\{N_i\}$ such that 
\eqnb
\eqna{y1}
\prb{y>y^{(N_i)}_C(t),\ \ \hat{y} ^{(N_i)} (y,t)-\hat{y}(y,t)> \eps_1} >\rho
\eqne
or
\[
\prb{y>y^{(N_i)}_C(t),\ \ \hat{y} ^{(N_i)} (y,t)-\hat{y}(y,t)< -\eps_1} >\rho
\]
holds.
We consider the case \eqnu{y1}, since the second case is dealt with 
similarly.
Let $y_1=\hat{y}(y,t)$.
With \prpu{ycytstrictincrease}, we have
\[ 
y=y_C(y_1,t) < y_C( y_1+\eps_1,t ).
\]
Let $\eps_2= y_C( y_1+\eps_1,t )-y >0$.
\prpu{trajectory} implies
\eqnb
\eqna{y3}
\prb{ |y_C^{(N_i)}(y_1+\eps_1,t )-y_C(y_1+\eps_1,t )|
 \geq \frac{\eps_2}{4}} \leq \frac{\rho}{2}
\eqne
holds for sufficiently large $N_i$.
Combining \eqnu{y1}, $y_1=\hat{y}(y,t)$ and \eqnu{y3}, 
we have
\[
\prb{ y>y_C^{(N_i)}(t),\ \ y_1+\eps_1 <\hat{y} ^{(N_i)} (y,t), \ \ 
|y_C^{(N_i)}(y_1+\eps_1,t )-y_C(y_1+\eps_1,t )|
< \frac{\eps_2}{4} } > \frac{\rho}{2}>0
\]
for sufficiently large $N_i$.
Note that \eqnu{HDLprf6} implies that 
if $y\le y'$ then 
$\dsp y\ssN_C(y,t)\le y\ssN_C(y',t)+\frac1N$\,.
This combined with \eqnu{yNChatyNC} implies  
that if $y>y_C^{(N_i)} (t)$ and $y_1+\eps_1 < \hat{y} ^{(N_i)} (y,t)$, then
\eqnb
\eqna{star1}
y_C ^{(N_i)} (y_1+\eps_1,t ) < y+\frac{\eps_2}{4}+\frac{3}{N_i} \,,
\eqne
for sufficiently large $N_i$.
On the other hand, 
if
 $\dsp |y_C^{(N_i)}(y_1 +\eps_1,t )-y_C(y_1 +\eps_1,t )| <\frac{\eps_2}{4} $,
 then 
\eqnb
\eqna{star2}
y+\eps _2 =y_C (y_1+\eps_1,t ) <
y_C^{(N_i)}(y_1 +\eps_1,t )+\frac{\eps_2}{4}\,.
\eqne
But \eqnu{star1} and \eqnu{star2} put together imply
\[y+\eps _2< y + \frac{\eps_2}{2}+\frac{3}{N_i}\,,\] 
which is a contradiction for large $N_i$.
Thus the assumption \eqnu{y1} is false,
which completes the proof of \thmu{HDL}.


\begin{thebibliography}{0}

\bibitem{longtail}
Chris Anderson, 
\textit{The Long Tail: 
Why the Future of Business Is Selling Less of More,}
 Hyperion Books, 2006.

\bibitem{Evans}
L.~C.~Evans, \textit{Partial differential equations,}
GSM \textbf{19}, AMS, 1998.

\bibitem{HH072}
K.~Hattori, T.~Hattori, 
\textit{Equation of motion for incompressible mixed fluid driven 
by evaporation and its application to online rankings,}
preprint, \url{http://arxiv.org/abs/0804.0330}

\bibitem{HH073}
K.~Hattori, T.~Hattori,
\textit{Mathematical analysis of long tail economy
using stochastic ranking processes,}
preprint, \url{http://arxiv.org/abs/0804.0321}

\bibitem{KL}
C.~Kipnis, C.~Landim,
\textit{Scaling limits of interacting particle systems,}
Springer, 1999.

\bibitem{OVY93}
S.~Olla, S.~R.~S.~Varadhan, H.~T.~Yau,
\textit{Hydrodynamical limit for a Hamiltonian system with weak noise,}
Commun.\ Math.\ Phys., \textbf{155} (1993) 523--560.

\bibitem{Spohn}
H.~Spohn,
\textit{Large scale dynamics of interacting particles,}
Springer, 1991.

\bibitem{VY97}
S.~R.~S.~Varadhan, H.~T.~Yau,
\textit{Diffusive limit of lattice gas with mixing conditions,}
Asian J.\ Math.\ \textbf{1} (1997) 623--678.


\end{thebibliography}
\end{document}